\documentclass{amsart}

\usepackage{latexsym} \usepackage{amssymb} \usepackage{amsfonts}
\usepackage{amscd} \usepackage{amssymb} \usepackage{amsmath}
\usepackage{url} 
\usepackage[all]{xy}

\newtheorem{prop}{Proposition}
\newtheorem{lem}[prop]{Lemma}
\newtheorem{cor}[prop]{Corollary} \newtheorem{thm}[prop]{Theorem}

\newtheorem{exmp}[prop]{Example}
\newtheorem{rem}[prop]{Remark} 

\def\OO{\mathcal{O}}
\def\lra{\longrightarrow }

\def\Xg{{\bf X}} 
\def\mm{\mathfrak{m}} 
\def\nn{\mathfrak{n}}

\def\Res{\mathrm{Res}} 

\def\ZZ{\mathbb{Z}}
\def\NN{\mathbb{N}}

\def\pp{\mathfrak{p}}
\def\Pp{\mathfrak{P}}

\def\SS{\mathfrak{S}}

\def\PP{{\mathbb{P}}}  
\def\NN{{\mathbb{N}}} \def\ZZ{{\mathbb{Z}}} 

\def\OO{{\mathcal{O}}}  

\def\Zc{{\mathcal{Z}}}  \def\Mc{{\mathcal{M}}}
 \def\Hc{{\mathcal{H}}}\def\Ic{{\mathcal{I}}}

 \def\Tg{\mathbf{T}} \def\Xg{\mathbf{X}}

\def\gr{\mathrm{gr}} \def\dim{\mathrm{dim}} 
 \def\deg{\mathrm{deg}} \def\Proj{\mathrm{Proj}}

\def\Syz{\mathrm{Syz}}
\def\codim{\mathrm{codim}}
\def\depth{\mathrm{depth}}
\def\Sym{\mathrm{Sym}} \def\ker{\mathrm{Ker}} 
 \def\Rees{\mathrm{Rees}}

\def\Res{{\mathrm{Res}}} 

 \def\mm{\mathfrak{m}} \def\pp{\mathfrak{p}}
 \def\indeg{\mathrm{indeg}}

\def\summ{\sum\nolimits}

\begin{document}

\title{Torsion of the symmetric algebra and implicitization}
\author{Laurent Bus\'e}
\address{Galaad, INRIA, 
2004 route des Lucioles, B.P. 93, 
06902 Sophia Antipolis Cedex, France}
\email{Laurent.Buse@inria.fr}

\author{Marc Chardin}
\address{Institut de Math\'ematiques de Jussieu, CNRS et Universit\'e Pierre et Marie Curie, 
4 place Jussieu, F-75252 Paris Cedex 05, France} 
\email{chardin@math.jussieu.fr}

\author{Jean-Pierre Jouanolou}
\address{Universit\'e Louis Pasteur,
  7 rue Ren\'e Descartes,
  67084 Strasbourg Cedex, France}
\email{jouanolo@math.u-strasbg.fr}

\date{\today}

\maketitle

\begin{abstract}
Recently, a method to compute the implicit equation of a pa\-ra\-me\-tri\-zed hypersurface 
has been developed by the authors. We address here some questions related to this method. First, we prove that the degree estimate  for the stabilization of the MacRae's invariant of $\Sym_{A}(I)_{\nu}$ is optimal. Then, 
we show that the extraneous factor that may appear in the process splits into a product a linear forms in the algebraic closure of the base field, each linear form being associated to a non complete intersection base point.
Finally, we make a link between this method and a resultant computation for the case of rational plane curves and  space surfaces. 
\end{abstract}

\section{Introduction}\label{intro}
Let $k$ be a field,  and $A:=k[X_1,\ldots,X_n]$ be a polynomial ring in $n\geq 2$ variables, with its standard $\NN$-grading  $\deg(X_i)=1$ for all
$i=1,\ldots,n$. Suppose given an integer $d\geq 1$ and $n+1$
homogeneous polynomials $f_1,\ldots,f_{n+1} \in A_d$, not all zero. We
will denote by $I$ the
ideal of $A$ generated by these polynomials and set $X:=\Proj (A/I)\subseteq \PP^{n-1}_k=\Proj(A)$. 
Let $T_1,\ldots,T_{n+1}$ be $n+1$ new indeterminates, $B:=k[T_1,\ldots,T_{n+1}]$ with its standard grading, and 
consider the rational map
\begin{equation}\label{lambda}
 \lambda :  \PP^{n-1}_k \rightarrow  \PP^{n}_k :  
 x  \mapsto  (f_1(x):f_2(x):\cdots:f_{n+1}(x)).
\end{equation}

We will always assume throughout this paper that $\lambda$ is generically finite onto its image, or equivalently that the closed image of $\lambda$ is a hypersurface $\Hc$, and that the base locus $X$ of $\lambda$ is supported on a finite set of points.
We will focus on the computation of an equation $H$ of the closed image $\Hc$  following a technique 
 developed in \cite{BuJo,BC03}. It is based on the computation of 
 MacRae's invariants of certain graded parts of the symmetric algebra
 $\Sym_A(I)$, which are $B$-modules, by means of the determinant of the corresponding graded
 parts of a  $\Zc$-approximation complex (see \cite{HSV} for the definition and a detailed study of these complexes). 
 
 The closure $\Gamma$ of the graph of $\lambda$ is a geometrically irreducible and reduced subscheme of 
 $\PP^{n-1}_k\times \PP^{n}_k$; its image under the projection $\pi : \PP^{n-1}_k\times \PP^{n}_k\rightarrow \PP^{n}_k$ is $\Hc$.
 
 The natural bigraded epimorphisms
 $$
 A\otimes_{k}B\buildrel{\phi}\over{\lra} \Sym_{A}\left(I(d)\right)\buildrel{\psi}\over{\lra} \Rees_{A}\left(I(d)\right)
 $$
 correspond to the inclusions of schemes
 $
 \Gamma \subseteq V \subseteq \PP^{n-1}_k\times \PP^{n}_k
 $
 where $V$ is the bi-projective scheme defined by $\Sym_{A}(I)$.
 
 By definition $\ker (\phi )$ is generated by the elements of bi-degree $(*,1)$ of the bigraded prime ideal $\Pp:=\ker (\psi\circ \phi )$. In other words, $\ker (\phi )$ is generated by the syzygies of the $f_{i}$'s : the elements $\sum_{i=1}^{n+1}a_{i}T_{i}\in A\otimes_{k}B_{1}$ such that $\sum_{i=1}^{n+1}a_{i}f_{i}=0$. 
 
As the construction of symmetric algebras and Rees algebras commute with localization, and both algebras are the quotient of a polynomial extension of the base ring by the Koszul syzygies on a minimal set of generators in the case of a complete intersection ideal, it follows that $\Gamma$ and $V$ coincide on $(\PP^{n-1}_k-S)\times \PP^{n}_k$, where $S$ is the (possibly empty) set of non locally complete intersection points of $X$. 
 
Recall that for a point $x:=\Proj(A/\pp )\in X$, the two notions of multiplicity :
\begin{itemize}
  \item the ``degree'', that we will denote by $d_{x}=d_{\pp}$,
which equals $\ell_{A_\pp} (A/I)_{\pp}=\ell_{\OO_{\PP^{n-1}_k,x}}(\OO_{X,x})$,
\item and the ``Hilbert-Samuel multiplicity'' $e_{x}=e_{\pp}$ given by the asymptotic behavior for $t\gg 0$ of $\ell_{A_\pp} (I^{t-1}/I^{t})_{\pp}=\ell_{\OO_{\PP^{n-1}_k,x}}(\Ic_{X,x}^{t-1}/\Ic_{X,x}^{t})$ (e.g.~\cite[\S 2.2]{BuJo}),
\end{itemize}
satisfies $e_\pp \geq d_\pp$ and this inequality is an
equality if and only if $I_\pp\subseteq A_{\pp}$ can be generated by a regular
sequence. More precisely, $e_\pp$ is the smallest degree of a $\pp$-primary complete intersection 
ideal contained in $I_\pp$, if the field $k$ is infinite  (see e.g.~\cite[\S 7]{Bou89}).

We  recall that if $M$ is a non zero $\ZZ$-graded $A$-module of finite type then its initial
degree, denoted $\indeg(M)$, is defined as the smallest integer $\nu$ such that the degree
$\nu$ part of $M$ is non-zero.

For $\nu\in \NN$, we define $\Sym_A(I)_\nu$ as the classes of elements of $B\otimes_{k}A_{\nu}$ in $\Sym_A(I)$. In other words, $\Sym_A(I)_\nu$ is the cokernel of the map :
$$
B\otimes_{k}\Syz (f_{1},\ldots ,f_{n+1})_{\nu +d}\lra B\otimes_{k}A_{\nu}
$$
sending $b\otimes (a_{1},\ldots ,a_{n+1})$ with $\deg (a_{i})=\nu$ (we follow here the usual convention for the grading of syzygy modules) to $b\otimes \sum_{i=1}^{n+1}a_{i}T_{i}$.

The following theorem is
proved in \cite{BuJo,BC03}, where the notation $\SS(-)$ stands for the MacRae's
invariant, as defined and studied in \cite{North}, of a  $B$-module. Recall that, as $B$ is regular, any finitely generated torsion $B$-module has a MacRae's invariant. If $M$ is a $B$-module that is not torsion, we set
 $\SS(M):=0$ (observe that the multiplicativity property of MacRae's invariants is preserved by this extension).

\begin{thm}[\cite{BuJo,BC03}]\label{TH}  Set $I_X:=I:_A\mm^\infty$. If $X$ is locally
  defined by at most $n$ equations\footnote{{\it i.e.}  $\Ic_{X,x}\subseteq \OO_{\PP^{n-1}_k,x}$ is generated by $n$ elements for any $x\in X$}  then,
  for any integer $\nu \geq \eta:=(n-1)(d-1)
  -\indeg(I_X),$ 
one has
  $$ \left( H^{\deg(\lambda)}\right) \supseteq \SS(\Sym_A(I)_\nu)=\SS(\Sym_A(I)_\eta)  \simeq
B\left(-d^{n-1}+\sum_{x\in X}d_{x}\right)$$
where the last isomorphism is a graded isomorphism of
$B$-modules. 
Moreover, 
the four following statements are equivalent
\begin{itemize}
    \item[\rm (i)] $I\subseteq A$ is of linear type on the punctured spectrum,
  \item[\rm (ii)] $X$ is locally a complete intersection,
   \item[\rm (iii)] $V=\Gamma$,
  \item[\rm (iv)] $ \left( H^{\deg(\lambda)} \right) = \SS(\Sym_A(I)_\nu)$ for
      all $\nu \geq \eta$.
  \end{itemize}
\end{thm}

Before going further, let us make two comments about the hypothesis of this theorem. First, it is proved in \cite{BC03} that if $X$ is not locally defined by $n$ equations at some base point $x$ (it is then locally defined by exactly $n+1$ equations) then the $\Zc$-approximation complex used to get a resolution of the symmetric algebra is no longer acyclic in high degrees (actually acyclicity and acyclicity in high degrees are equivalent in this case). Moreover, notice that in this case, $\pi (V)=\PP^{n}_k$ as $V$ contains $\{ x\} \times \PP^{n}_k$. Second, the degree bound $\eta$ is sharp in the following sense

\begin{prop}\label{optdeg} Under the hypotheses of Theorem \ref{TH},
$$\SS(\Sym_A(I)_\nu)=\SS(\Sym_A(I)_\eta) \text{ if and only if } \nu \geq \eta.$$
\end{prop}

We will use the following lemma in the proof,
\begin{lem}\label{JI} If  $J\subset I$ is a complete intersection ideal generated
by $n-1$ forms of degree $d$ and $l\in \ZZ$, then
$$
l <\indeg (I_X)\ \Leftrightarrow (I_X)_{l}=0\ \Leftrightarrow \ (I_X/J)_{l}=0.
$$
\end{lem}

\begin{proof} [Proof of Lemma \ref{JI}]  Notice that as $J$ is unmixed and $\dim (A/J)>0$, it suffices to show that $J_{l}\not= (I_X)_{l}$ for $l =\indeg (I_X)$. This is clearly the case if $\indeg (I_X)<d=\indeg (J)$. If $\indeg (I_X)=d$, $J_{d}= (I_X)_{d}$ implies that $J_d=I_d$,  which in turn shows that $\lambda$ factors through a rational map to $\PP^{n-2}_{k}$, contradicting the assumption
that the closed image of $\lambda$ has dimension $n-1$.
\end{proof}

\begin{proof} [Proof of Proposition \ref{optdeg}] The case $n=2$ follows immediately from \cite[\S 2]{BC03}, so we can assume that $n\geq 3$. By faithfully flat base change $k\rightarrow \overline{k}$, we may assume that $k$ is algebraically closed, and in particular infinite.  As  $\Sym_A(I)_\nu =0$ for $\nu <0$ it suffices to prove that if $\nu \geq 0$ then $$\SS(\Sym_A(I)_\nu)=\SS(\Sym_A(I)_\eta) \Leftrightarrow \nu \geq \eta:=(n-1)(d-1)-\indeg(I_X).$$

Let $\nu\geq 0$,  $\mu :=(n-1)(d-1)-\nu$ and $D:=\deg (\det (\Zc_{\eta}))(=d^{n-1}-\sum_{x\in X}d_{x})$.
As $\codim (I)=n-1$, $I$ is generated in degree $d$ and $k$ is infinite, there exists a complete intersection ideal $J\subset I$ 
generated by $n-1$ forms of degree $d$. We will now prove that, if  $\Zc_{\nu}$ is generically acyclic, one has
$
\deg (\det (\Zc_{\nu}))=D-\dim (I_X/J)_{\mu -1}.
$
This together with Lemma \ref{JI} establishes our claim.

Recall that for a finitely generated graded module $M$ one has $H_{M}(l)=P_{M}(l)+\sum_{i}(-1)^{i}h^{i}_{\mm}(M)_{l}$ ($l\in \ZZ$), where $P_M$, respectively $H_M$, denotes the Hilbert polynomial, respectively the Hilbert function, of $M$ and $h^{i}_{\mm}(M)_{l}:=\dim_{k}H^{i}_{\mm}(M)_{l}$. We will denote by $Z_i$ and
$H_i$, the $i$-th module of cycles and $i$-th module of homology  of the Koszul complex $K(f_1,\ldots,f_{n+1};A)$, respectively. The degree $\nu$ part of the $\Zc$-complex  given by the $f_i$'s will be denoted by $\Zc_{\nu}$. Recall that $(\Zc_i)_{\nu}=B\otimes_k (Z_i)_{\nu +id}$.

As $\mu \leq (n-1)(d-1)\leq (n-1)d$, \cite[4.1]{BC03} shows that $H_{Z_{0}}(\nu ) = P_{Z_{0}}(\nu )$,
\begin{align} \nonumber
H_{Z_{1}}(\nu +d) &= P_{Z_{1}}(\nu +d)+h^{0}_{\mm}(H_{1})_{\mu +d-1}-H_{H_{2}}(\mu +d-1),\\ \nonumber
H_{Z_{2}}(\nu +2d) &= P_{Z_{2}}(\nu +2d)+h^{0}_{\mm}(H_{0})_{\mu -1}-H_{H_{1}}(\mu -1)+H_{H_{2}}(\mu -1),\\ \nonumber
H_{Z_{i}}(\nu +id) &= P_{Z_{i}}(\nu +id)+\summ_{j}(-1)^{i+j}H_{H_{j}}(\mu -(i-2)d-1) \ \text{ for } \  3\leq i\leq n 
\end{align}
(notice that $H_p$ coincides with $Z_p$ in degree at most $(p+1)d-1$, and therefore $H^{n}_{\mm}(Z_{i})_{\nu +id}\simeq
(Z_{n-i})_{\mu -(i-2)d-1}^{*}$ coincides with $(H_{n-i})_{\mu -(i-2)d-1}^{*}$).

Set $\chi (l):=\sum_{j}(-1)^{j}H_{H_{j}}(l)$ ($l\in \ZZ$) and remark that 
\begin{align} \nonumber
H_{Z_{2}}(\nu +2d) &= P_{Z_{2}}(\nu +2d)+\chi (\mu -1)-H_{R/I_X}(\mu -1),\\ \nonumber
H_{Z_{i}}(\nu +id) &= P_{Z_{i}}(\nu +id)+(-1)^{i}\chi (\mu -(i-2)d-1) \ \text{ for } \ 3\leq i\leq n.
\end{align}
Now $\sum_{i} (-1)^{i}P_{Z_{i}}(l+id)=0$ for any $l$ and $P(l):=\sum_{i} (-1)^{i+1}iP_{Z_{i}}(l+id)$  is constant 
and equal to $D$. If $\Zc_{\nu}$ is generically acyclic, then $ \sum_{i} (-1)^{i}H_{Z_{i}}(\nu +id)=0$ 
and we can write, using our expressions for $H_{Z_{i}}$,
\begin{align} \nonumber
\deg (\det (\Zc_{\nu}))&=\sum_{i} (-1)^{i-1}iH_{Z_{i}}(\nu +id)\\ \nonumber
&=P(\nu )+\sum_{i} (-1)^{i-1}(i-1)[ H_{Z_{i}}(\nu +id)-P_{Z_{i}}(\nu +id)]\\ \nonumber
&=D+H_{R/I_X}(\mu -1)-\sum_{i\geq 0}(i+1)\chi (\mu -id-1).
\end{align}
 It follows that 
$\deg (\det (\Zc_{\nu}))=D+H_{R/I_X}(\mu -1)-\sum_{i\geq 0}(i+1)\chi (\mu -id-1).$

The Hilbert series of $\chi$ only depends on $d$ (the degree of the $f_{i}$'s) and $n$ and is given by
$$
\sum \chi (l)t^{l}=\frac{(1-t^d)^{n+1}}{(1-t)^{n}},
$$
which shows that the series corresponding to $F(l):=\sum_{i\geq 0}(i+1)\chi (l -id)$ is
$$
\frac{(1-t^d)^{n+1}}{(1-t)^{n}}\sum_{i\geq 0}(i+1)t^{id}=\frac{(1-t^d)^{n-1}}{(1-t)^{n}}.
$$
Hence $F=H_{R/J}$ and $\deg (\det (\Zc_{\nu}))=D+H_{R/I_X}(\mu -1)-H_{R/J}(\mu -1)$.
\end{proof} 

\begin{rem} Under the hypotheses of Theorem \ref{TH}, set  $\varepsilon :=\indeg(I_X)$. One can expect that $\Sym_A(I)_{\eta -1}$ is not a torsion $B$-module unless $I$ is $n$-generated. This is indeed the case if $I$ is saturated, or if $I$ has no syzygy of degree
less then $d+\varepsilon$, or if $I_X=A$ (i.e. $X=\emptyset$). 
\end{rem}

It is known that the degree of the hypersurface $\Hc$, closure of the image of $\lambda$, is given by the formula (see for instance \cite[Section 6]{SUV} for an algebraic proof, or  \cite[2.5]{BuJo} for a more geometrical approach) 
$$\deg(\lambda)\deg(\Hc)=d^{n-1} - \sum_{x \in X}{e_{x}}.$$
Therefore, under the assumption of Theorem \ref{TH}, for all $\nu\geq \eta$ the MacRae's invariant
$\SS(\Sym_A(I)_\nu)$ is a principal ideal of $B$ whose generator
(unique up to multiplication by a non-zero constant in $k$) equals
$H^{\deg(\lambda)}$ times an \emph{extraneous homogeneous factor} $G\in B$ whose total degree is equal to 
$\sum_{x \in X} (e_x-d_x )$
and whose support is the codimension 1 part of $\pi (V-\Gamma )$.

\medskip

In the following section \ref{exfct} we give an explicit description of this extraneous
factor $G$. Then, in section \ref{rescp} we relate the computation of the MacRae's invariant $\SS(\Sym_A(I)_\nu)$ to a resultant computation in the case  
of algebraic plane curves ($n=2$) and space surfaces ($n=3$). 

\section{Description of the extraneous factor}\label{exfct}

Recall that $S\subset X$ is the set of points where $X$ is not locally a complete intersection. Assume that $x=\Proj(A/\pp )\in S$ is a base point  locally defined by exactly $n$ equations. Then one of the generators of $I$ may be written in terms of the other ones in any sufficiently small open neighborhood of $x$ (or equivalently in $A_{\pp}$), and this unique relation (up to multiplication by a function invertible in $A_{\pp}$) specializes at $x$, thus defining a non zero linear form $L_{x}\in B_{1}\otimes_k k(x)$. It turns out that the extraneous factor $G$ can be described in terms of these linear forms.

\begin{prop}\label{propG} If the field $k$ is algebraically closed and $X$ is locally generated by at most $n$ elements then, up to multiplication by a non-zero constant in $k$, then 
  $$G=\prod_{x \in S} L_x^{e_x-d_x} \in B.$$
Moreover, all syzygies $a_{1}T_{1}+\cdots +a_{n+1}T_{n+1}\in A\otimes_{k}B_{1}\simeq A^{n+1}$ specializes at $x$ to a multiple of $L_{x}$, and at least one
  does not specialize to $0$.
\end{prop}
\begin{proof} Let $K$ be the kernel of the canonical epimorphism $\Sym_A (I) \rightarrow \Rees_A(I)$. The multiplicativity of MacRae's invariants shows the equality of ideals of $B$
$$ \SS(\Sym_A(I)_\nu)= \SS(\Rees_A(I)_\nu)\SS(K_\nu)$$
for all integers $\nu \geq 0$. Since $\SS(\Rees_A(I)_\nu))=(H^{\deg(\lambda)})$ for all $\nu\geq 0$ (see for instance the proof of \cite[Theorem 2.5]{BuJo}) and $ \SS(\Sym_A(I)_\nu)=(GH^{\deg(\lambda)})$ for all $\nu \geq \eta$ by Theorem \ref{TH}, we deduce that $\SS(K_\nu)=(G)$ for all $\nu \geq \eta$. 

Fix an integer $\nu \geq \eta$. The divisor $\SS(K_\nu)=(G)$ is supported on $\pi (V-\Gamma )$ and we know that $V$ and $\Gamma$ coincides outside $S\times \PP^{n}_k$. Moreover, $V$ is defined by linear forms in the $T_{i}$'s, in particular $V-\Gamma$ is supported on a union of linear spaces over the points of $S$, therefore its closure is a scheme whose support is a union of linear spaces $\{ x\} \times V_{x}$ for $x\in S$. Now  $V_{x} \subseteq W_{x}:=\{ L_{x}=0\}$. It follows that $G=\prod_{x \in S} L_x^{\mu_x}$ for some integers $\mu_x \geq 0$. 
The following lemma applied to $A$ localized at the prime corresponding to $\{ x\} \times W_{x}$ shows that $\mu_x=e_x-d_x$.
%
\end{proof}

\begin{lem} Let $(R,\nn )$ be a Cohen-Macaulay local ring of dimension $d$, $I$ be a $\nn$-primary ideal of $R$ and $K$ be the kernel of the canonical map $\Sym_R (I) \rightarrow \Rees_R (I)$. If $I$ is generated by $d+1$ elements and the $\Zc$-complex of $I$ is acyclic, then $K$ has dimension $\leq d$ and its multiplicity in dimension $d$ is equal to  
$e_I(R)-l_R(R/I)$.
\end{lem}
\begin{proof}
The commutative exact diagram 
$$\xymatrix{
&0\ar[d]&0\ar[d]&&\\
&K(1)\ar^{\lambda_K}[r]\ar[d]&K\ar[d]&&\\
&\Sym_R (I)(1)_{\geq 1}\ar^{\lambda}[r]\ar[d]&\Sym_R (I)\ar[r]\ar[d]&\Sym_{R/I} (I/I^2)\ar[r]\ar[d]&0\\
0\ar[r]&\Rees_R (I)(1)_{\geq 1}\ar^{can}[r]\ar[d]&\Rees_R (I)\ar[r]\ar[d]&\gr_{I} (R)\ar[r]\ar[d]&0\\
&0&0&0&\\}
$$
shows that $\ker (\lambda_K )\simeq \ker (\lambda )$ and gives rise to an exact sequence of 
graded $\Sym_R (I)$-modules
$$\xymatrix{
0\ar[r]&\ker (\lambda )\ar[r]&K(1)\ar^{\lambda_K}[r]&K\ar[r]&\Sym_{R/I} (I/I^2)\ar[r]&\gr_{I} (R)\ar[r]&0.\\
}$$
Furthermore, if the $\Zc$-complex of $I$ is acyclic, then $H_1({\Mc })\simeq \ker (\lambda )$ and $H_i ({\Mc })=0$ for $i\geq 2$. Hence the above exact sequence of graded
$\Sym_R (I)$-modules of finite type supported on $V(I)$ implies an equality of Poincar\'e series :
$$
(1-t^{-1})S(K,t)=S(H_0({\Mc}),t)-S(H_1({\Mc}),t)-S(\gr_{I} (R),t)
$$
in  the $\ZZ[[t]][t^{-1}]$-module $K_0(R/I)[[t]][t^{-1}]$.
If $I$ is generated by $a:=(a_1,\ldots ,a_m)$, one has 
$$
\sum_i (-1)^i S(H_i({\Mc}),t)
={\frac{1}{(1-t)^m}}\sum_{i}(-1)^i S(H_i(a;R),t)t^i.
$$
This gives 
$$
S(K,t)={\frac{t}{(1-t)}}S(\gr_{I} (R),t)-{\frac{t}{(1-t)^{m+1}}}\sum_{i}(-1)^i S(H_i(a;R),t)t^i.
$$

Set $d:=\dim R$. As $I$ is $\nn$-primary, $K_0(R/I)=\ZZ$ and $S(M,t):=\sum_{n\in \ZZ}l_R(M_n)t^n$ for a
finitely generated graded $R/I$-module $M$. Recall that $S(\gr_{I} (R),t)={\frac{Q(t)}{(1-t)^d}}$ with
$Q(1)=e_I (R)$ and $S(K,t)={\frac{P(t)}{(1-t)^{d+1}}}$ for some $P\in \ZZ [t]$. It follows that
$$
(*)\quad\quad (1-t)^{m-d}P(t)=t(1-t)^{m-d}Q(t)-t(\sum_{i=0}^{m}(-1)^i l_R(H_i(a;R))t^i).
$$
If further $R$ is Cohen-Macaulay and $m=d+1$, then $H_j(a;R)=0$ for $j\not\in \{ 0,1\}$. Setting
$t=1$ in formula $(*)$ shows that $l_R(H_1(a;R))=l_R(H_0(a;R))=l_R(R/I)$, hence
$$
P(t)=tQ(t)-{\frac{t}{1-t}}(l_R(H_0(a;R))-l_R(H_1(a;R))t)=t(Q(t)-l_R(R/I))
$$
which gives $P(1)=e_I(R)- l_R(R/I)$.
\end{proof}

\section{Link with some resultants}\label{rescp}

We begin this section by focusing on the
particular case where the ideal $I\subset A=k[X_1,\ldots,X_{n}]$ generated by the polynomials
$f_1,\ldots,f_{n+1}$ (that define the map $\lambda$ in \eqref{lambda}) has projective dimension 1. From Hilbert-Burch Theorem (see for instance \cite[20.15]{Eis}), $I$ admits a finite free graded resolution of the form 
\begin{equation}\label{resI}
0\rightarrow \bigoplus_{i=1}^{n} A(-\mu_i-d) \xrightarrow{ M } A(-d)^{n+1}
\xrightarrow{\left[ f_1 \ f_2 \ \cdots \ f_{n+1} \right]} A \rightarrow A/I
\rightarrow 0
\end{equation}
where $M=\left[p_{i,j} \right]_{\substack{ i=1,\ldots,n+1, \ j=1,\ldots,n}}$, with the properties that $\sum_{i=1}^{n}\mu_i=d$ and the ideal of $n$-minors of $M$ is equal to $I$. 

 When the ideal $I$ has projective dimension 1 our approach to the implicitization
problem via Theorem \ref{TH} coincides with a resultant
computation.
As a consequence of the Auslander-Buchsbaum formula, the ideal $I$ may have projective dimension 1 only if $n=2$ or $n=3$ (remember that we assumed that $\dim(A/I)\leq 1$).  

More precisely, define the polynomials 
$$L_i:=\sum_{j=1}^{n+1} p_{j,i}T_j \in A[T_1,\ldots,T_{n+1}] \ \text{ for all } \ i=1,\ldots,n$$
that are respectively bi-homogeneous of bi-degree $(\mu_i;1)$ in $A[\Tg]=k[\Xg][\Tg]$. We denote by $\Res(L_1,\ldots,L_{n})$ the resultant of $L_1,\ldots,L_n$ w.r.t.~the homogeneous variables $X_1,\ldots,X_n$. 

\begin{prop}[$n=2,3$]\label{resdet} Assume that the ideal $I$ is locally generated by at most $n$ elements and that \eqref{resI} is a free resolution of $I$. Then, for all $\nu\geq \eta$  
  $$(\Res(L_1,\ldots,L_{n}))=\SS(\Sym_A(I)_\nu) \subset B.$$
\end{prop}
\begin{proof} Let $J$ be the ideal generated by polynomials
  $L_1,\ldots,L_{n}$ in the polynomial ring $A[T_1,\ldots,T_{n+1}]$ and denote by $Q$
  the quotient algebra $A[T_1,\ldots,T_{n+1}]/J$. 

On the one hand,  by definition of the symmetric algebra 
  $Q\simeq \Sym_A(I)$ as bi-graded modules.  On the other hand, 
 $\Res(L_1,\ldots,L_{n})=\SS(Q_\nu)$, for all $\nu \gg 0$, if the sequence
 $(L_1,\ldots,L_n)$ is regular in $A[T_1,\ldots,T_{n+1}]$ (see e.g.~\cite[\S 3.5]{J95}). 
  We deduce that the claimed result is proved if we show that $\Sym_A(I)\simeq Q$ is a complete
 intersection. Recall Avramov's criterion \cite{Avramov}: 
{\it $\Sym_A(I)$ is a complete intersection if and only if, for all integer $r\in \{1,\ldots,n\}$, $\depth_{I_r(M)}(A)\geq n-r+1$, where $ I_r(M)$ denotes the ideal of $A$ generated by all the $r\times r$-minors of the matrix $M$} (notice that the fact that $B$ is a complete intersection also follows from \cite[2.1]{ASV81}).

In the case $n=2$, $M$ is a $3\times 2$ matrix and, by Hilbert-Burch Theorem, $\depth_{I_2(M)}(A)\geq 2$. Since $I_2(M) \subset I_1(M)$,  $\depth_{I_1(M)}(A)\geq \depth_{I_2(M)}(A)$ and we obtain that $\Sym_A(I)$ is a complete intersection in this case.

 For $n=3$, the Hilbert-Burch Theorem shows that $\depth_{I_3(M)}(A)=\depth_{I}(A) \geq 2$ since we assumed that $I$ defines isolated points. It follows that 
$$\depth_{I_1(M)}\geq
\depth_{I_2(M)}\geq \depth_{I_3(M)} \geq 2.$$
Moreover, as $I$ is supposed to be locally generated by at most 3 elements, basic properties of Fitting ideals \cite[\S 20.2]{Eis} imply that $V(I_1(M))=\emptyset$, i.e.~$\depth_{I_1(M)}(A) \geq 3$. Therefore, $\Sym_A(I)$ is also a complete intersection in this case\footnote{Notice that if $I$ is locally minimally generated by four elements at some point, then $V(I_1(M))$ is non-empty, $\depth_{I_1(M)}(A)=2$ and hence $\Res(L_1,\ldots,L_n)= 0$ in $B$.}. 
\end{proof}

In the case $n=2$ ($\lambda$ parametrizes a plane
curve), the ideal $I$ always admits a free resolution of the
form \eqref{resI} and is always locally generated by at most 2 elements; therefore,
Proposition \ref{resdet} shows that as soon as $\lambda$ is generically finite, the
resultant of $L_1$ and $L_2$ always compute an implicit equation
of the image of $\lambda$ to the power $\deg(\lambda)$, without any
extraneous factor. 

In the case $n=3$ ($\lambda$ parametrizes a space 
surface), the situation is a little more intricate. The ideal $I$
admits a resolution of the form \eqref{resI} if it defines some isolated
base points (at least one) and if it is saturated with respect to the maximal ideal
$(X_1,X_2,X_3)$ of $A$ (see for instance \cite[\S 20.4]{Eis}). Proposition \ref{propG} and Proposition \ref{resdet} prove 
the following result that has been conjectured in 
\cite[Conjecture 5.1]{BCD}:

\begin{cor}\label{coro} If the field $k$ is algebraically closed, the ideal $I$ has codimension 2, is saturated and is locally generated by at most 3 elements, then $I$ admits a resolution as in \eqref{resI} and, with the above notation, one has\footnote{Observe that this equality is rigorously an equality of divisors since $H$ and the $L_\pp$'s are defined up to multiplication by a nonzero element in $k$ and have not been explicitly chosen.}
$$\Res(L_1,L_2,L_3)= H(T_1,\ldots,T_4)^{\deg(\lambda)}\prod_{x \in V(I)} L_x (T_1,\ldots,T_4)^{e_x -d_x},$$
where $H \in B$ denotes an implicit equation of the image of  $\lambda$.
\end{cor}

\begin{exmp} Consider the following $4\times 3$-matrix 
$$M:=\left[\begin{array}{ccc}
X_1^2+X_3^2 & X_2 & X_2X_3 \\
X_1^2-X_1X_2 & X_1+2X_2 & X_1^2+X_3^2 \\
X_1X_2+X_1X_3 & 2X_2+X_3 & X_2^2 \\
X_1X_3+X_2^2 & X_3 & X_1X_2
\end{array}\right]$$
whose ideal $I_3(M) \subset A$  of $3\times 3$ minors is of codimension 2. This matrix fits into  \eqref{resI} and hence defines the parametrization of a surface $\Hc$ in $\PP^3$.
Since $V(I_1(M))=\emptyset$, Corollary \ref{coro} applies and the resultant of $L_1,L_2,L_3$ is  the product of two irreducible polynomials:  
an implicit equation of $\Hc$ which has degree 6, and an extraneous factor $G:=T_2^2 + 2T_3T_4 + T_3^2 + T_4^2$. 
As predicted by Proposition \ref{propG}, $G$ splits into a product of two linear forms over the complex numbers field $\mathbb{C}$:
$$G=(T_2-iT_3-iT_4)(T_2+iT_3+iT_4) \in \mathbb{C}[T_1,T_2,T_3,T_4]$$
where $i\in \mathbb{C}$,  $i^2=-1$. Also, evaluating the matrix $M$ at the two projective points $(X_1:X_2:X_3)=(1:0:\epsilon i)$, $\epsilon=\pm 1$, we find 
$$M(1:0:\epsilon i):=\left[\begin{array}{ccc}
0 & 0 & 0 \\
1 & 1 & 0 \\
\epsilon i & \epsilon i & 0 \\
\epsilon i & \epsilon i & 0
\end{array}\right].$$
Hence these two points are base points of the parametrization (the last column is the null vector) which are almost complete intersection since the rank of $M(1:0:\epsilon i)$ is 1. The linear factors of  $G$  can be read from a non-zero column of $M(1:0:\epsilon i)$. 
\end{exmp}

From now on we will bring the case $n=3$ into focus. It is then possible to weaken the hypotheses of Proposition \ref{resdet} using the fact that the first syzygy module of $I=(f_1,f_2,f_3,f_4)$ becomes a locally free module of rank 3 after localization by any non-zero linear form $l\in A_1$ (by Hilbert's syzygy theorem applied to localizations of $A_{(l)}\simeq k[X_1,X_2 ]$), hence a free $A_{(l)}$-module of rank 3, by Quillen-Suslin Theorem. For simplicity, assume that $l=X_3$. Then $I_{(X_3)} \subset A_{(X_3)}\simeq k[X_1,X_2]$ admits a resolution of the form
\begin{equation}\label{resIl}
0\rightarrow  {A_{(X_3)}}^3 \xrightarrow{ M } {A_{(X_3)}}^{4}
\xrightarrow{\left[ \tilde{f}_1 \ \tilde{f}_2 \ \tilde{f}_3 \ \tilde{f}_{4} \right]} I_{(X_3)}
\rightarrow 0
\end{equation}
where  $\tilde{f}_i(X_1,X_2)=f_i(X_1,X_2,1)$ (observe that $\max_i(\deg(\tilde{f}_i))=d$, since we assumed that $I$ has codimension at least 2) and $M=[p_{i,j}]_{i=1,\ldots,4, \, j=1,\ldots,3}$. 
For all $j=1,2,3$, we define $L_j(X_1,X_2;\Tg)=\sum_{i=1}^4 p_{i,j}T_i \in A_{(X_3)}[X_1,X_2]$ and its homogenization in the variable $X_3$ w.r.t.~the variables $X_1,X_2$ 
$$L_j^h(\Xg;\Tg)=\sum_{i=1}^4 X_3^{\deg_{X_1,X_2}(L_j)-\deg_{X_1,X_2}(p_{i,j})}p_{i,j}T_i=\sum_{i=1}^4p_{i,j}^h(\Xg)T_i \in A[\Tg].$$
Set $d_j=\deg_\Xg(L_j^h)=\deg_{X_1,X_2}(L_j)$ and $M^h=\left[ p^h_{i,j} \right]_{i=1,\ldots,4, \, j=1,\ldots,3}.$
%

\begin{prop}\label{affres} Assume $I$ is locally generated by at most $3$ elements off $V(X_3)$. With the above notation, consider the two conditions:
\begin{itemize}
\item[(a)] $V(I_1(M^h),X_3)=\emptyset$ and $V(X_3) \nsubseteq V(I_2(M^h))$,
\item[(b)] there exists homogeneous polynomials $P_1,P_2,P_3 \in k[X_1,X_2]$ without common factor and a  linear form $Q \in k[\Tg]$ such that  $L_i^h(X_1,X_2,0,\Tg)=P_iQ$ for all $i=1,2,3$. 
\end{itemize}
Then
\begin{align*}
\Res(L_1^h,L_2^h,L_3^h)\neq 0 \in B \iff 
\text{\rm (a) holds or (b) holds.} 
\end{align*}
Furthermore, if  $k$ is algebraically closed  and 
$\Res(L_1^h,L_2^h,L_3^h)$ is nonzero then it is equal to 
 $$H(\Tg)^{\deg(\lambda)}  \prod_{x\in V(I_2(M^h))\setminus V(X_3)} L_x(\Tg)^{e_x-d_x} \prod_{x\in V(I_2(M^h))\cap
 V(X_3)} l_x(\Tg)^{\mu_x}$$
where $l_x$ is a linear form such that $\langle l_x\rangle_k=\langle L_1^h(x),L_2^h(x),L_3^h(x)\rangle_k$, and 
$$\sum_{x\in V(I_2(M^h))\cap
 V(X_3)}\mu_x = d_1d_2+d_2d_3+d_1d_3-d^2+\sum_{x \in V(I)\setminus V(X_3)}d_x.$$
\end{prop}

\begin{proof} We begin by proving the characterization for $\Res(L_1^h,L_2^h,L_3^h)\neq 0$. Observe that if  (b) holds then  $V(X_3)\subseteq V(I_2(M^h))$; in particular (a) and (b) do not overlap. Also Avramov's criterion implies that $(L_1^h,L_2^h,L_3^h)$ is a $A[\Tg]$-regular sequence if and only if  (a) holds, in which case $\Res(L_1^h,L_2^h,L_3^h)\neq 0$.

 ($\Leftarrow$): We may assume that (b) holds. For $i=1,2,3$, 
$$L_i^h(\Xg,\Tg)=P_i(X_1,X_2)Q(\Tg) + X_3R_i(\Xg,\Tg).$$
As $(L_1,L_2,L_3)$ form a $A_{(X_3)}[\Tg]$-regular sequence, $Q\neq 0$ and $V(P_1,P_2,P_3)=\emptyset \subset \PP^1$ we deduce that $\Res(L_1^h,L_2^h,L_3^h)\neq 0$.

 ($\Rightarrow$): We proceed by contradiction and assume that  both (a) and (b) do not hold. If $x\ \in V(I_1(M^h),X_3)\neq \emptyset$ then  $L_i^h(x;\mathbf{T})=0$, $i=1,2,3$, hence $\Res(L_1^h,L_2^h,L_3^h)=0$. Else, since the specialization $X_3=0$ makes $M^h$ a matrix of rank $\leq 1$ and the ring $ k[X_1,X_2][\Tg ]$ is factorial, there exists polynomials $(P_i(X_1,X_2))_{i=1,2,3}$ and $Q(X_1,X_2,\Tg)$ such that 
\begin{align*}
L_i^h(X_1,X_2,0,\Tg) =P_i(X_1,X_2)Q(X_1,X_2,\Tg) .
\end{align*}
Since (b) does not hold, either the $P_i$'s have a common root, or  $Q$ has a positive degree w.r.t.~the variables $X_1,X_2$. It follows that $\Res(L_1^h,L_2^h,L_3^h)= 0$ since $L_i^h \in (P_1,P_2,P_3,X_3)\cap (Q,X_3)$, $i=1,2,3$.

We now turn to the proof of the second statement. Let $\nu\geq \eta$ be an integer. From resultant theory it follows that $\Res(L_1^h,L_2^h,L_3^h)$ is a generator of the principal $B$-ideal $\SS(A[\Tg]_\nu/(L_1^h,L_2^h,L_3^h)_\nu)$. Consider the canonical exact sequence of $B$-modules
\begin{multline}\label{SSplit}
 0 \rightarrow H^0_{(X_3)}(A[\Tg]_\nu/(L_1^h,L_2^h,L_3^h))_\nu  \rightarrow A[\Tg]_\nu/(L_1^h,L_2^h,L_3^h)_\nu \xrightarrow{\phi} \\ \Sym_A(I)_\nu/H^0_{(X_3)}(\Sym_A(I))_\nu \rightarrow 0.
\end{multline}
From Proposition \ref{propG} and its proof, we deduce that the MacRae's invariant of the last term of \eqref{SSplit} is generated by
$$H(\Tg)^{\deg(\lambda)}  \prod_{x\in V(I_2(M^h))\setminus V(X_3)} L_x(\Tg)^{e_x-d_x}.$$
Furthermore, the MacRae's invariant of the first term of \eqref{SSplit} is supported above each point $x \in V(X_3)\subset \PP^2$ such that the linear space $\{ L_1^h(x)=L_2^h(x)=L_3^h(x)=0\}$ has codimension one (see for instance \cite{HR}); we deduce that it is generated by $\prod_{x\in V(I_2(M^h))\cap
 V(X_3)} l_x(\Tg)^{\mu_x}$. This proves the claimed formula for $\Res(L_1^h,L_2^h,L_3^h)$ by multiplicativity of MacRae's invariants applied on \eqref{SSplit}. The last claim  is a straightforward computation.
\end{proof}

The above proposition answers questions raised in \cite[\S5]{CCL05}. Let us emphasize that $\Res(L_1^h,L_2^h,L_3^h)$ and a generator of $\SS(\Sym_A(I)_\nu)$, with $\nu \geq \eta$, may only differ by linear forms which are associated to base points supported on the line $V(X_3)$. 

\begin{exmp} Consider the example \cite[Example 2.2]{CCL05} of a canonical Steiner surface.
The first syzygy module of the $\tilde{f}_i$'s is given by the matrix
$$M=\left[
\begin{array}{ccc}
 0 & 0 & 1\\
X_1X_2 & 1+X_2^2 & -X_1 \\
1+X_1^2 & X_1X_2 & 0 \\
-2X_1 & -2X_2 & 0
\end{array}
\right].$$
One can verify that $L_1^h,L_2^h,L_3^h$ is a regular sequence by Avramov's criterion,  and compute
$$\Res(L_1^h,L_2^h,L_3^h)={T_2}^4\, H(T_1,T_2,T_3,T_4),$$
where $H=0$ is an implicit equation of the Steiner surface. 
Observe that the extraneous factor $T_2^4$ is associated to the base point $X_1=X_3=0$. Also, notice that $\SS(\Sym_A(I)_\nu)$ is generated by $H$ for all $\nu\geq \eta=2$.
\end{exmp}

The choice of the matrix $M$ is crucial: from a given $M$ it is easy to get another choice such that $\Res(L_1^h,L_2^h,L_3^h)=0$. 
This leads to the following question: is it always possible to find a matrix $M$ such that $\Res(L_1^h,L_2^h,L_3^h)\neq 0$ ? The answer likely uses the condition that $\sum_{i=1}^3\deg_\Xg(L_i)$ has to be minimal.


\begin{thebibliography}{1}

\bibitem{ASV81}
J.F.~Andrade, A.~Simis, and W.~Vasconcelos.
\newblock On the grade of some ideals.
\newblock {\em Manuscripta Math.}, 41:241--254, 1981.


\bibitem{Avramov}
Luchezar~L. Avramov.
\newblock Complete intersections and symmetric algebras.
\newblock {\em J. Algebra}, 73(1):248--263, 1981.

\bibitem{Bou89}
Nicolas Bourbaki.
\newblock {\em \'{E}l\'ements de math\'ematique}.
\newblock Masson, Paris, 1983.
\newblock Alg\`ebre commutative. Chapitre 8. Dimension. Chapitre 9. Anneaux
  locaux noeth\'eriens complets. [Commutative algebra. Chapter 8. Dimension.
  Chapter 9. Complete Noetherian local rings].



\bibitem{BC03}
Laurent Bus{\'e} and Marc Chardin.
\newblock Implicitizing rational hypersurfaces using approximation complexes.
\newblock {\em J. Symbolic Comput.}, 40:1150-1168, 2005.

\bibitem{BCD}
Laurent Bus{\'e}, David Cox, and Carlos D'Andrea.
\newblock Implicitization of surfaces in {${\mathbb P}\sp 3$} in the presence
  of base points.
\newblock {\em J. Algebra Appl.}, 2(2):189--214, 2003.

\bibitem{BuJo}
Laurent Bus{\'e} and Jean-Pierre Jouanolou.
\newblock On the closed image of a rational map and the implicitization
  problem.
\newblock {\em J. Algebra}, 265(1):312--357, 2003.

\bibitem{CCL05}
Falai Chen, David Cox, and Yang Liu.
\newblock The {$\mu$}-basis and implicitization of a rational parametric surface.
\newblock {\em  J. Symbolic Comput.}, 39(6):689--706, 2005.

\bibitem{Eis}
David Eisenbud.
\newblock {\em Commutative algebra}, volume 150 of {\em Graduate Texts in
  Mathematics}.
\newblock Springer-Verlag, New York, 1995.
\newblock With a view toward algebraic geometry.


\bibitem{HSV}
J. Herzog, A. Simis and W. Vasconcelos. 
\newblock Koszul homology and blowing-up rings.  
\newblock {\em Commutative algebra (Trento, 1981)},  pp. 79--169, Lecture Notes in Pure and Appl. Math., 84, Dekker, New York, 1983.


\bibitem{HR}
C. Huneke and M. Rossi.
\newblock The dimension and components of symmetric algebras.
\newblock {\em J. Algebra}, 98:200--210, 1986.

\bibitem{J95}
Jean-Pierre Jouanolou.
\newblock Aspects invariants de l'\'elimination.
\newblock {\em Adv. Math.}, 114(1):1--174, 1995.

\bibitem{North}
D.~G. Northcott.
\newblock {\em Finite free resolutions}.
\newblock Cambridge University Press, Cambridge, 1976.
\newblock Cambridge Tracts in Mathematics, No. 71.

\bibitem{SUV}
A. Simis, B. Ulrich and W. Vasconcelos. 
\newblock Codimension, multiplicity and integral extensions.  
\newblock {\em Math. Proc. Cambridge Philos. Soc.},  130  (2001),  no. 2, 237--257.

\end{thebibliography}

\end{document}